\documentclass[11pt,twoside]{article}

\usepackage[total={6.5in, 9in}, top=1in, left=1in]{geometry}
\usepackage{graphicx}
\usepackage{graphics}
\usepackage{comment}
\usepackage{amsmath}
\usepackage{amsfonts}
\usepackage{amssymb}
\usepackage{euscript}
\usepackage{amsthm}
\usepackage{hyperref}
\usepackage{enumerate}
\usepackage{tikz}
\usetikzlibrary{decorations}
\usetikzlibrary{decorations.pathmorphing}
\usetikzlibrary{decorations.pathreplacing}
\usetikzlibrary{decorations.markings}
\usetikzlibrary{decorations.text}
\usetikzlibrary{arrows.meta}
\usetikzlibrary{calc}
\usetikzlibrary[arrows,positioning, shapes, fit]
\usetikzlibrary{snakes}
\usepackage{tkz-euclide}
\newtheorem{thm}{Theorem}[section]
\newtheorem{dfn}[thm]{Definition}
\newtheorem{prop}[thm]{Proposition}
\newtheorem{cor}[thm]{Corollary}
\newtheorem{lem}[thm]{Lemma}

\newtheorem{obs}[thm]{Observation}
\newtheorem{alg}[thm]{Algorithm}

\newtheorem{open}[thm]{Open Question}
\newcommand{\Z}{\operatorname{Z}}

\definecolor{grn1}{rgb}{0,0.9,0}
\definecolor{grn2}{rgb}{0,0.65,0}
\definecolor{grn3}{rgb}{0,0.35,0}

\definecolor{turq}{rgb}{0,0.5,.5}

\begin{document} 
\title{Well-forced graphs}
\author{Cheryl Grood\thanks{Swarthmore College}, Ruth Haas\thanks{Smith College and University of Hawaii, Manoa}, Bonnie Jacob\thanks{National Technical Institute for the Deaf, Rochester Institute of Technology (bcjntm@rit.edu)}, Erika King\thanks{Hobart and William Smith Colleges}, and Shahla Nasserasr\thanks{College of Science, Rochester Institute of Technology}}
\maketitle
\begin{abstract}
A graph in which all  minimal zero forcing sets are in fact minimum size is called {\em well-forced}.
This paper characterizes well-forced trees and presents an algorithm for determining which trees are well-forced. Additionally, we characterize which vertices in a tree are contained in no minimal zero forcing set.  
\end{abstract}


\section{Introduction}\label{Sec:Intro}

 Zero forcing is a topic that has received much recent attention in the mathematical literature. In this paper, we add to the body of knowledge about zero forcing in trees in two distinct ways. We characterize trees with all minimal zero forcing sets the same size; and we characterize vertices in trees that occur in no minimal zero forcing set.
 
 Zero forcing was introduced in \cite{MR2388646} in 2008 as a combinatorial technique for finding an upper bound on the maximum multiplicity of the eigenvalues of symmetric matrices whose nonzero off-diagonal entries represent the edges of the graph with no restriction on the diagonal entries. Since then the zero forcing parameter and its variants 
have been studied 
extensively; see \cite{HLS2022inverse} for more details. Zero forcing can be seen as a game with a given ``color change rule,"  where the goal of a player is to color all vertices  blue. In this paper, we consider standard zero forcing which involves only two colors and one player. The standard color change rule (CCR) is as follows: if a vertex $v$ is  blue and each of its neighbors except the neighbor $u$ is blue, then the vertex $u$ is changed to blue. The game starts with some initial subset of vertices colored blue  and the rest white.  At any step, if a vertex can be colored blue  by the CCR, then the player applies the rule, otherwise a new vertex must be added to the initial set. The player typically aims for the minimum number of vertices that need to be initially colored blue in order to color  all of the vertices of the graph blue. An initial set of  blue vertices that can result in all of the vertices of the graph being colored blue is called a {\em zero forcing set}. A {\em minimum zero forcing set} is a zero forcing set of minimum cardinality, and the {\em zero forcing number} of a graph is the cardinality of a minimum zero forcing set in the graph. The notation $\Z(G)$ is used for the zero forcing number of a graph $G$. It is known that a graph may have several minimum zero forcing sets. For instance either of the degree-one vertices of a path is a zero forcing set. A {\em minimal zero forcing set} is a zero forcing set such that the  removal of any vertex from the set results in a set that is no longer a zero forcing set.  We note two facts about minimal zero forcing sets: first, there may be multiple minimal zero forcing sets in one graph, second, they are not necessarily minimum zero forcing sets. For example, any two adjacent vertices of degree two on a path $P_n; n\geq 5$ form a minimal zero forcing set. In this work, we study graphs for which every minimal zero forcing set is also minimum. We refer to these graphs as {\em well-forced graphs}. 

Extensive research has been done characterizing graphs with the property that all sets of some defined type that are minimal are also minimum. A graph is called  well-$X$ for a parameter $X$  if  every minimal  $X$ set  is also minimum, or in other words, if every minimal set has the same cardinality. (Note that the same statements hold if we replace ``minimal'' and ``minimum'' with ``maximal'' and ``maximum'' respectively).
For example, well-covered graphs are graphs in which every maximal independent set is maximum \cite{MR1677797, MR1254158}; well-dominated graphs are graphs in which every minimal dominating set is minimum  \cite{MR4197372, MR0942485}; and well-indumatchable graphs are graphs in which all maximal induced matchings are maximum \cite{MR4453102}. Questions that are NP-Complete or co-NP-Complete for certain graph parameters often become solvable in polynomial time on these ``well" graphs, making them important classes to understand. Here our definition of well-forced graphs is an analog of this research.


The main focus of this paper is to study well-forced graphs. A recent study of a  related problem is in \cite{BorisCarlson}, where the authors provide several families of graphs that have exponentially many minimal zero forcing sets, and consider how a graph can have minimal zero forcing sets of different sizes. They do not address  the question of a structural description of well-forced graphs. We address that question in this paper.  In Section \ref{Sec:Intro}, we introduce the problem and  provide basic terminology used in this paper. Critical to our understanding of  well-forced graphs is a characterization of  irrelevant vertices, which are those that appear in no minimal zero forcing sets.     In Section \ref{Sec:Gratis ch} we build up results on paths, generalized stars and finally trees to give a characterization of irrelevant vertices in trees. 
Finally, in Section \ref{Sec:well trees}, we  begin the study of well-forced graphs, culminating in a characterization of  well-forced trees. In addition, we provide an algorithm that determines if a tree is well-forced.

\subsection{Preliminaries}
 All graphs in this paper are finite and  simple.  That is, $G=(V,E)$ where $V$ is a finite set of vertices and there is at most one (undirected) edge between each  pair of distinct vertices. Vertices $u, v\in V$ are {\em adjacent} or {\em neighbors} if the edge $uv \in E$. The {\em degree} of a vertex is the number of neighbors it has. The set of neighbors of $v$ is denoted by $N(v)$. A {\em pendent vertex} or \emph{pendant} is a vertex with degree one. 
For other standard graph theory terminology the reader may consult  \cite{CZ}. We proceed with some specialized definitions.

\begin{dfn}
We say that pendent vertices $u, w$ in a graph $G$ are \emph{double pendants} if $N(u)=N(w)=\{v\}$. We say that $v$ \emph{has a double pendant}.  
\end{dfn}

\begin{dfn}A vertex $v \in V(G)$ is of \emph{high degree} in $G$ if $\deg_G(v) \geq 3$.\end{dfn}

A \emph{generalized star} is a tree with at most one vertex of high degree.  A \emph{pendent path} is an induced subgraph of a graph $G$ that is a path with a vertex of degree one in $G$, and any other vertices on the path have degree two in $G$.  



\begin{dfn}
A \emph{pendent generalized star} of a graph $G$ is an induced subgraph $R$ of $G$ such that
\begin{enumerate}
\item there is exactly one vertex $v$ of $R$ that is of high degree in $G$ 
\item $G - v$ has $k+1$ components and exactly $k$ of the components 
 of $G-v$ are pendent paths of $v$.
\item $R$ is induced by the vertices of the $k$ pendent paths and $v$.
\end{enumerate}
\end{dfn}

In a generalized star $T$, the \emph{length} of a leg is the number of edges on the path from the center vertex of $T$ to the pendent vertex on the leg. The same statement holds for a pendent generalized star $R$.  

\begin{figure}[htbp]
\begin{center}
\begin{tikzpicture}[auto, scale=0.8]
\tikzstyle{vertex}=[draw, circle, inner sep=0.9mm]
\node (v1) at (0,0) [vertex, fill=grn3]{};
\node (v2) at (0.7,-0.7) [vertex]{};
\node (v3) at (-0.7,-0.7) [vertex]{};
\node (right1) at (1,0) [vertex]{};
\node (right2) at  (2,0) [vertex, fill=grn2]{};
\node (rightlower) at (2.7, -0.7) [vertex]{};
\node (rightupper1) at (2.7, 0.7) [vertex]{};
\node (rightupper2) at (3.4, 1.4) [vertex, fill=grn1]{};
\node (rightupper3) at (4.1, 2.1) [vertex]{};
\node (left1) at (-1, 0) [vertex]{};
\node (left2) at (-2, 0) [vertex, fill=grn1]{};
\node (left2lower) at (-2, -1) [vertex, fill=grn3]{};
\node (leftdouble1) at (-2.7, -1.7) [vertex]{};
\node (leftdouble2) at (-1.3, -1.7) [vertex]{};
\node (veryleft1) at (-3, 0) [vertex]{};
\node (veryleft2) at (-4, 0) [vertex, fill=grn2]{};
\node (veryleft3) at (-5, 0) [vertex]{};
\node (veryleftup1) at (-4, 1) [vertex]{};
\node (veryleftup2) at (-4, 2) [vertex]{};
\node (leftdiagup1) at (-4.7, 0.7) [vertex]{};
\node (leftdiagup2) at (-5.4, 1.4) [vertex, fill=grn3]{};
\node (leftdiagup3) at (-6.1, 2.1) [vertex]{};
\node (leftdiagdouble1) at (-5.4, 2.4) [vertex]{};
\node (leftdiagdouble2) at (-6.4, 1.4) [vertex]{};

\draw (v2) -- (v1) -- (v3);
\draw (v1) -- (right1) -- (right2) -- (rightlower);
\draw (right2) -- (rightupper1) -- (rightupper2) -- (rightupper3);
\draw (v1) -- (left1) -- (left2) --(left2lower);
\draw (leftdouble1) -- (left2lower) -- (leftdouble2);
\draw (left2) -- (veryleft1) -- (veryleft2) -- (veryleft3);
\draw (veryleft2) -- (veryleftup1) -- (veryleftup2);
\draw (veryleft2) -- (leftdiagup1) -- (leftdiagup2) -- (leftdiagup3);
\draw (leftdiagdouble1) -- (leftdiagup2) -- (leftdiagdouble2);
\end{tikzpicture}
\end{center}
\caption{A tree with  $B_0$, $B_1$, and $B_2$ vertices shaded \textcolor{grn3}{dark}, \textcolor{grn2}{medium}, and \textcolor{grn1}{light} green respectively.}
\label{figure:bvertices}
\end{figure}
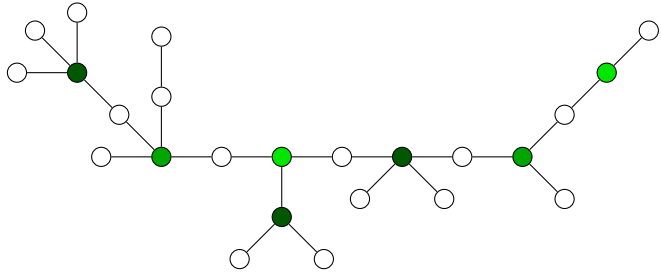

We now define $B$-vertices.  
\begin{dfn}
In a graph $G$, if $v \in V(G)$ has a double pendant, then we say that $v \in B_0$.   Remove from $G$ all vertices in $B_0$ and their pendent neighbors to create the graph $G_1$. Then $B_1$ is the set of vertices in $G_1$ that each have a double pendant.   Continuing, for $i \geq 1$, we construct $G_{i+1}$ by removing from $G_i$ the set $B_i$ and all neighbors of $B_i$ that have degree one in $G_i$.  The vertices in $B_0, B_1, B_2, \ldots$ are said to be \emph{$B$-vertices} of $G$.  The pendent neighbors that are removed along with $B_i$ are called the  \emph{pseudoleaves} of $B_i$.  If $u \in B_i$ and $v$ is a pendent neighbor in $G_i$ of $u$, then we say that $v$ is a \emph{pseudoleaf of} $u$. 
\end{dfn}
A tree with $B$-vertices shaded is shown in Figure \ref{figure:bvertices}.

\begin{dfn}
A vertex $v \in V(G)$ that has the property that $v$ is in no minimal zero forcing set of $G$ is an \emph{irrelevant} vertex. 
\end{dfn}
 In this paper we will show that in a tree, a vertex is irrelevant if and only if it is a $B$-vertex, and that the reverse implication is true for general 
 graphs. Irrelevant vertices were first defined in \cite{BongCarlsonCurtis},  where they were so named
 because they are irrelevant in the 
 construction of the reconfiguration graph for zero forcing. 

\subsection{Basic properties of zero forcing sets}
A zero forcing set $S$ has an associated ordered list of forces. 
A \emph{forcing chain}  is a sequence of vertices, $(v_1, v_2, \ldots, v_k)$ such that $v_i$ forces $v_{i+1}$ for each $i$, $1 \leq i \leq k-1$.  Note that $k=1$ is permitted.  An ordered list of forces under $S$ produces a set of forcing chains whose union is all vertices of $G$.
A \emph{maximal forcing chain} is a forcing chain that is not a proper subset of another chain using the same 
list of forces.  This implies that $v_1 \in S$ and that $v_k$ does not force any other vertex under this  ordered 
list of forces.  Throughout this paper, when we refer to forcing chains, we always assume that they are maximal.  

It is important to note that one zero forcing set may have more than one associated ordered list of forces and also
different sets of forcing chains.  We call one particular set of forcing chains under a zero forcing set $S$ a \emph{realization} of forcing chains. 

 A \emph{path cover} of a tree $T$ is a set of vertex-disjoint paths that are induced subgraphs of $T$ such that each vertex of $T$ is in some path of the path cover. A \emph{minimum path cover}, denoted by $P(T)$, of a tree $T$ is a path cover with the smallest possible number of paths among all path covers of $T$.   We use the following facts from \cite{MR2388646} throughout the paper: by taking an end vertex from every path in a minimum path cover of a tree, we produce a minimum zero forcing set; also, the set of forcing chains produced from a minimum zero forcing set in a tree is a minimum path cover. We make use of the following related lemma throughout this paper.

\begin{lem}
Given a path cover of a tree, let $S$ be a set of vertices consisting of either an end vertex or a pair of adjacent internal vertices from each path in the path cover. Then $S$ is a zero forcing set.   \label{lem:pathcovers}
\end{lem}

\begin{proof}
We proceed by induction on the number of paths in the path cover, similarly to the proof of \cite[Proposition 4.2]{MR2388646} but without the restriction that the path cover be minimum.  If the path cover contains just one path, then $T$ is a path, and the result holds.

Suppose that for up to $k$ paths in the path cover, the result holds.  Consider a tree $T$ with a path cover consisting of $k+1$ paths.  Because $T$ is a tree, there exists a path in the path cover that is only joined to the rest of the tree by one edge not on the path.  Choose such a path and call it $P$.  Let $uv$ be the edge that joins $P$ to the rest of the tree, $T'$, with $v \in V(P)$.

 Take either end vertex of $P$, or a pair of adjacent internal vertices in $P$ (if $P$ is long enough), and do the same for all other paths in the path cover. Call this set $S$.  Either $v \in S$, or $v$ will be forced by a neighbor on $P$. Thus the remaining tree $T'$ does not depend on $P$, and the induction hypothesis gives us that all vertices in $V(T')$ are forced, including $u$.  Any remaining vertices on $P$ will be forced, giving us that $S$ is a zero forcing set of $T$.  
\end{proof}

Given a zero forcing set $S$ and a set of forcing chains associated with the set, the \emph{reversal} of $S$ is the set consisting of the last vertex in each forcing chain.  In \cite{barioli2010zero}, the reversal was shown to be a zero forcing set; since the reversal has the same cardinality as the original zero forcing set, the reversal of a minimum zero forcing set is also a minimum zero forcing set.  The reversal of a minimal zero forcing set is not necessarily a minimal zero forcing set, however.  

\section{Characterization of irrelevant vertices}\label{Sec:Gratis ch}
In this section, we present results that culminate in a characterization of irrelevant vertices in trees.  Most results presented in this section apply only to trees, though some apply to graphs in general.  First, we note a simple condition that guarantees that a vertex is irrelevant.

\begin{prop} 
Suppose $u,w \in V(G)$ form a double pendant with $N(u)=N(w)=\{v\}$.  Then under any minimal zero forcing set $S$, there is a realization of forcing chains such  that one of $\{u,w\}$ forces $v$.  Consequently, at least one of $\{u,w\}$ is in $S$, and $v \notin S$. \label{prop:pendent_not_minimal}
\end{prop}

\begin{proof}
For any minimal zero forcing set $S$, at least one of $u,w$ must be in $S$; otherwise neither will be forced.   Without loss of generality, suppose $u \in S$.  Then $u$ forces $v$; hence $v \notin S$.  
\end{proof}

\subsection{Irrelevant vertices in paths and generalized stars}

Next we establish some basic properties of minimal zero forcing sets and irrelevant vertices in paths and generalized stars.  
\begin{lem} \label{lem:patheveryvertexminimal} 
If $n \neq 3$, then every vertex in $P_n$ is in a minimal zero forcing set.
\end{lem}

\begin{proof}
Either end vertex of the path forms a minimum zero forcing set. Let $v$ be an internal vertex in $P_n, n \geq 4$.  Then $v$ has at least one neighbor $u$ that is also an internal vertex, and $\{u,v\}$ forms a minimal zero forcing set. 
\end{proof}

If $S$ is a zero forcing set of a generalized star $T$, then all legs but at most one must contain a vertex from $S$.  If all legs contain a vertex from $S$, then $S$ is not minimal, giving us the following observation.  

\begin{obs}
If $T$ is a generalized star with $\ell$ legs and minimal zero forcing set $S$, then exactly $\ell-1$ legs of $T$ contain a vertex from $S$.  
\label{obs:allbutoneleg}
\end{obs}

We now present a condition under which a generalized star has no irrelevant vertices.  

\begin{lem} \label{lem:genstareveryvertexminimal}
If $T$ is a generalized star with at most one leg of length one, then every vertex in $T$ is in some minimal zero forcing set. \end{lem}
\begin{proof}
By Lemma \ref{lem:patheveryvertexminimal} we know that the result holds if $T$ is a path, so we can assume $T$ has exactly one vertex $c$ with $\deg(c) \geq 3$.  Let $\ell=\deg(c)$.
By Observation \ref{obs:allbutoneleg}, for any minimal zero forcing set $S$ of $T$, exactly $\ell - 1$ legs of $T$ have a vertex from $S$.  Also note that $\Z(T)=\ell -1$.  

Let $v \in V(T)$.  If $v$ is a leaf, take $S$ to consist of $v$ and $\ell-2$ other leaves.  Then $S$ is a minimum and hence minimal zero forcing set.  

If $v=c$, take $S$ to be $v$ and $\ell-1$ non-leaf neighbors of $v$.  
Then $v$ forces the leg with no vertices from $S$, and each other vertex in $S$ forces its respective leg; we cannot remove $v$ from $S$ because then each vertex in $S$ has two white neighbors, and we cannot remove any other vertex from $S$ because then  two legs of $T$ contain no vertex from $S$.  Hence $S$ is a minimal zero forcing set.  This construction also applies if $v$ is a non-leaf neighbor of $c$. 

The only remaining case is that $v$ is a non-leaf vertex at distance at least two from $c$.  In this case, construct $S$ as follows: add to $S$ the vertex $v$ and a neighbor $u$ of $v$ that is not a leaf.  On each of  $\ell-2$ remaining legs that are not of length one, add the vertex that is adjacent to $c$.  Since $u, v$  force the full leg they are on and $c$, we see that $S$ is a zero forcing set.  We show that $S$ is minimal: if we remove $u$ or $v$, no color changes occur because every vertex in $S$ has two white neighbors.  If we remove another vertex from $S$, then  two legs have no vertices from $S$.  Thus, $S$ is a minimal zero forcing set.
\end{proof}

\subsection{Pendent generalized stars and irrelevant vertices}
We now establish results related to pendent generalized stars that will be helpful in characterizing irrelevant vertices in trees. The next lemma is similar to Observation \ref{obs:allbutoneleg}, but applies to pendent generalized stars.

\begin{lem}\label{lem:starleguse} 
Let $G$ be a graph with pendent generalized star $R$, where $R$  has $\ell$ legs. 
Any  zero forcing set of $G$ must contain a  vertex from at least $\ell-1$ legs of $R$.
\end{lem}

\begin{proof}
Suppose $\ell_1, \ell_2$ are two legs of $R$ with no vertex from set $S$. 
 The only neighbor of a vertex on  $l_i$ that is not on $l_i$ is the center vertex $c$, but $c$ can never force $l_1$ and $l_2$ because it has  at least two white neighbors, one on each of $l_1, l_2$.  So $S$ is not a zero forcing set.
\end{proof}

The next two lemmas and corollary focus on the case that the center of a pendent generalized star is in a minimal zero forcing set.  
\begin{lem} \label{centerminimalhasneighbors}
Let $G$ be a graph with pendent generalized star $R$, where $R$ has $\ell$  legs and center vertex $c$.  If $c \in S$ for some minimal zero forcing set $S$, then $V(R) \cap S$ consists of $c$ and either $\ell-1$ or $\ell$ neighbors of $c$, all of whom have degree two. \label{lem:neighborsofcenter}
\end{lem}

\begin{proof}
Since $c\in S$, no leaves of $R$ are in $S$, else minimality is contradicted.
Consequently, $R$ has at most one leg of length one, since $S$ has at least one vertex from at least $\ell-1$ legs of $R$ by Lemma \ref{lem:starleguse}, and in a length-one leg, the only vertex is the leaf. 

Suppose  $v \in S$ where $v$ is a degree-two vertex in $R$ that is not adjacent to $c$.  Let $u$ be a neighbor of $v$.  If $u, v \in S$, then $u,v$ force the leg they're on, and will also force $c$, contradicting minimality of $S$ since then we could remove $c$ from $S$ and still produce a zero forcing set.  Thus, no neighbor of $v$ is in $S$.   If $u$ is forced by $v$, then the other neighbor of $v$ must be forced at some point or in $S$, which means that it could force $v$, contradicting minimality of $S$ since we could remove $v$.  Similarly, if $u$ is forced by its other neighbor (if $u$ is not a leaf), then $u$ can force $v$, contradicting minimality of $S$.  

The only remaining vertices, then, on $R$ that can be in $S$ are the degree-two neighbors of $c$.  Since $S$ must have a vertex from at least $\ell-1$  legs of $R$  by Lemma \ref{lem:starleguse}, and the degree-two neighbors of $c$ are our only candidates, the proof is complete.
\end{proof}

By replacing a center vertex  and its neighbor on a pendent generalized star with the leaf on the same leg in Lemma \ref{centerminimalhasneighbors} to produce a smaller zero forcing set, we note the following corollary of Lemma \ref{lem:neighborsofcenter}.
\begin{cor}
Let $G$ be a graph with $R$ a pendent generalized star and $c$ the center of $R$.  If $c \in S$ for some zero forcing set $S$, then $S$ is not a minimum zero forcing set.  \label{cor:nocenter}
\end{cor}

\begin{lem} \label{lem:centerminimallegsminimal}
Let $G$ be a graph with  pendent generalized star $R$. Let $c$ be the center of $R$ and suppose $c \in S'$ for some minimal zero forcing set $S'$ of $G$.  Then for all $v \in V(R)$, there exists a minimal zero forcing set $S$ such that $v \in S$. 
\end{lem}

\begin{proof}
Since $c \in S'$,  the pendent generalized star $R$ has at most one leg of length one by Proposition \ref{prop:pendent_not_minimal}, and by Lemma \ref{centerminimalhasneighbors}, we have that  $S' \cap V(R)$ must consist of $c$ and either $\ell -1$ or $\ell$ neighbors of $c$ on $R$, where $\deg(c)=\ell+1$.

Note that we can construct a minimum zero forcing set containing any leaf  in any pendent generalized star, since if a minimum zero forcing set does not contain a leaf, then its reversal does.  Since $S'$ can contain any non-leaf neighbors of $c$, we need only show the result for vertices in $R$ that have degree two and are at least distance two from $c$.  Let $v$ be such a vertex, and let $u$ be a neighbor of $v$ of degree two.  Let $w$ be the vertex that is adjacent to $c$ on the leg containing $v$.  (Note it is possible that $u=w$).

If $w \in S'$, define $S= \left( S' \backslash \{c, w\} \right) \cup \{u, v\}$.  Since $u, v$ force all vertices on their leg and $c$ as well, and $S'$ is a zero forcing set, we see $S$ is also a zero forcing set.  For minimality, note that if we remove $u$ or $v$, then $c$ will not be forced. (By Lemma  \ref{centerminimalhasneighbors}, all other legs of $R$ contain at most one vertex in $S'$ and hence in $S$.   Moreover, if some vertex in $T \backslash R$ could force $c$, then this would contradict the minimality of $S'$).  Also by minimality of $S'$, $S \backslash \{x\}$ is not a zero forcing set for any $x \in S' \cap S$.  Hence, $S$ is a minimal zero forcing set containing $v$.

If $w \notin S'$, let $w'$ be a degree-two neighbor of $c$ in $R$ with $w' \in S'$.  Note we can swap $w$ with $w'$ in $S'$ and still produce a minimal zero forcing set.  We can proceed then with the construction in the case $w \in S'$ to complete the proof. 
\end{proof}

We now show that under certain conditions, we can remove legs from a pendent generalized star in a tree, and find that any minimal zero forcing set on the remaining subtree is a subset of a minimal zero forcing set on the original tree.  

\begin{lem} \label{lem:tprime}
Suppose $T$ is a tree with pendent generalized star $R$ where $R$ has at most one leg of length one.  Let $T'$ be the subtree of $T$ with all legs of $R$ removed except for a shortest leg.  Then for any minimal zero forcing set $S'$ of $T'$, there exists a minimal zero forcing set $S$ of $T$ with $S \cap V(T') = S'$.
\end{lem}

\begin{proof}
In $T$, label the legs of $R$ from $1$ to $\ell$ where Leg $\ell$ is a shortest leg.  Denote the respective leaves  by $u_1, u_2, \ldots, u_{\ell}$. Denote the respective neighbors of $c$ by $w_1, w_2, \ldots, w_{\ell}$.  Note that $u_{\ell}=w_{\ell}$ is possible if $R$ has a leg of length one, but for $i < \ell$, $u_i \neq w_i$.  
 
Let $S'$ be a minimal zero forcing set of $T'$. Note that $\deg_{T'}c=2$.  Let $P$ denote the path from $u_{\ell}$ to $w_{\ell}$.  Note that $T'$ contains a pendent path from $u_{\ell}$ through $x$, the neighbor of $c$ not on $R$; the path $P$ is a subpath of this pendent path, and it's possible that $|V(P)|=1$.

In $T$, let $S_1 = S' \cup \{w_1, w_2, \ldots, w_{\ell-1}\}$.  Let $S_2 = S_1 \backslash \{w_1\}$.   Note that if $\ell=2$, then $S_2=S'$.  Also note that  $S_2 \cap V(T') = S_1 \cap V(T') = S'$.  We will show that one of these two sets $S_1$, $S_2$ is a minimal zero forcing set of $T$, which will complete the proof.

 We have either: (1) $v \in S'$ for some $v \in V(P)$, or (2) $v \notin S'$ for any $v \in V(P)$.

\underline{Case 1:} Suppose $v \in S'$ for some $v \in V(P)$.  There are three possibilities: (a) $v= u_{\ell}$; (b) $u,v \in S'$ for some $u \in V(P)$ with $u, v \neq u_{\ell}$ where $uv \in E(T)$; or (c) $v=w_{\ell} \neq u_{\ell}$.

First consider (a).  In $T'$, we see that $c \notin S'$ since that would contradict minimality as $u_{\ell}$ will force $P$ and then $c$.  Also, $c$ is not forced by its neighbor $x$ not on $P$, since then $c$ would force $P$, and we need not have $u_{\ell} \in S'$.  

We claim that $S_1$ or $S_2$ is a minimal zero forcing set on $T$.  First, note that $S_1$ is a zero forcing set of $T$.  If $S_2$ is a zero forcing set, then it must be minimal by Lemma \ref{lem:starleguse} since only $\ell-1$ legs of $R$ have a vertex from $S_2$, and if we remove a vertex $y$ from $S_2 \cap V(T')$ and still produce a zero forcing set of $T$, then we could remove $y$ from $S'$ and produce a zero forcing set of $T'$, contradicting minimality.  If $S_2$ is not a zero forcing set, then we show that $S_1$ must be a minimal zero forcing set.  We know we cannot remove $w_i$ from $S_1$, $1 \leq i \leq \ell-1$, because $S_2$ is not a zero forcing set.  If we remove $u_{\ell}$, then $c$ will not be forced because we know $x$ doesn't force $c$ under $S'$, and $S_1\cap V(T') =S'$.  Finally, if we can remove any other vertex $y$ from $V(T') \cap S_1$ and still produce a zero forcing set on $T$, then $S' \backslash \{y\}$ is a zero forcing set on $T'$ contradicting minimality of $S'$.

For (b), the same construction and arguments hold.

For (c), note first that the neighbor of $w_{\ell}$ on $P$ is not in $S'$ because that case was covered in (b).  Somehow, $c$ must be forced under $S'$.  If  $x$ forces $c$, then $c$ could force $v=w_{\ell}$ violating minimality, so we must have that $c \in S'$. Also note that if the neighbor $x$ of $c$ is in $S'$ or forced by some vertex other than $c$, then $c$ could force $w_{\ell}$.  Hence, $c$ must be required to force $x$ under $S'$ in $T'$.  We claim that $S_1$ is a minimal zero forcing set on $T$. $S_1$ is clearly a zero forcing set.  For minimality, note that $c$ must force its neighbor $x$, and therefore can't force any neighbors on $R$.  Thus, removing $w_i$ for any $i$ results in a set that is not a zero forcing set.  If we remove $c$, since we know that $c$ is required to force $x$, no color changes on $R$ or to $x$ will happen.  Finally, by minimality of $S'$, we cannot remove any vertex in $V(T')\cap S'$ from $S_1$, completing the proof of Case 1.

\underline{Case 2:} Suppose no vertex from $P$ is in $S'$.  Consider $S_1$ on $T$.  We can see that $S_1$ is a zero forcing set, since $S'$ will force all of $T'$ including $c$, allowing $w_i$ to force Leg $i$ for each $i <\ell$.   The set $S_1$ is also minimal.  Only $\ell-1$ legs of $R$ have a vertex from $S_1$, so none of these vertices may be removed by Lemma \ref{lem:starleguse}.  If we can remove a vertex from $S_1\cap V(T')$ and produce a zero forcing set of $T$, then we can remove the same vertex from $S'$ to produce a zero forcing set of $S'$, contradicting minimality of $S'$ and completing the proof.   
\end{proof}

\subsection{Irrelevant vertices in trees}
We now use the above results to establish through the following theorem  that if a tree has an irrelevant vertex, then it must have 
double pendants. 

\begin{thm} \label{thm:nodoublependants}
 Let $T$ be a tree with no double pendants.  Then every vertex in $T$ is in a minimal zero forcing set of $T$.
 \end{thm}

\begin{proof}
 We proceed by induction on $k$, the number of high-degree vertices in $T$.  For $k=0$ and $k=1$, the statement holds by Lemmas \ref{lem:patheveryvertexminimal} and \ref{lem:genstareveryvertexminimal} respectively.

Assume that if $T$ has up to $k$ high-degree vertices and no double pendants, then every vertex in $T$ is in some minimal zero forcing set of $T$.  

Let $T$ be a tree with $k+1$ high-degree vertices and no double pendants where $k \geq 1$.  By \cite{fallat2007minimum},  the tree $T$ must have a pendent generalized star $R$ with center $c$.  
Construct $T'$ from $T$ as in the statement of Lemma \ref{lem:tprime} using $R$.  The tree $T'$ has $k$ high-degree vertices (since $\deg_{T}c \geq 3$  but $\deg_{T'}c =2$).  Note that $c$ has at most one pendent neighbor in $T'$ because $T$ has at least two high-degree vertices, and there are no double pendants in $T'$.
Thus, applying the induction hypothesis to $T'$,  every vertex in $T'$ is in some minimal zero forcing set $S'$ of $T'$.  By Lemma \ref{lem:tprime}, then, every vertex of $T'$ is in some minimal zero forcing set $S$ of $T$ as well.  It remains to show that all vertices in $T \setminus T'$ are in some minimal zero forcing set of $T$.  
By Lemma \ref{lem:centerminimallegsminimal}, since $c$ is in some minimal zero forcing set of $T$, each vertex in $R$ is in some minimal zero forcing set of $T$, completing the proof. 
\end{proof}

We note the following observation that applies to disconnected graphs that we use implicitly throughout.
\begin{obs}
The set $S$ is a minimal zero forcing set of a graph $G$ if and only if  $S \cap V(H)$ is a minimal zero forcing set of $H$ for each connected component $H$ of $G$.  \label{obs:disconnected}
\end{obs}

We next establish that removal of a vertex with a double pendant along with all pendent neighbors of the vertex leaves minimal zero forcing sets largely unchanged. 
\begin{thm}
  Let $T$ be a tree with vertex $v$ that has leaf neighbors $L = \{l_1, \dots , l_k\}$ where $k\geq2$. Define $T' = T-\{v, l_1, \dots , l_k\}$. $S$ is a minimal zero forcing set of $T$ if and only if one of the following holds. \label{thm:starremovalstillminimal2}
\begin{enumerate}
\item $S$ contains exactly $k-1$ vertices of $L$, and $S'=S \cap T'$ is a minimal zero forcing set of $T'$. \label{case:kminus1}
\item $S=\left( S' \backslash \{u\} \right)  \cup L $  where $uv \in E(T')$, and $S'$ is a minimal zero forcing set of $T'$ with $u \in S'$.  \label{case:k}
\end{enumerate}
\end{thm}
\begin{proof}
We begin with the reverse direction.  First, assume that Statement \ref{case:kminus1}
holds.  That is, assume that $S'$ is a minimal zero forcing set of $T'$.  Then for any $j$, $1 \leq j \leq k$, under $S=S'\cup L - \{l_j\}$ in $T$, a vertex in $L\backslash\{l_j\}$ forces $v$, and $v$ forces $l_j$, and since $S'$ is a zero forcing set of $T'$, we have that $S$ is a zero forcing set of $T$.  For minimality of $S$, note that we cannot remove $l_i$ from $S$ for any $i \neq j$, since then $v$ will have two leaf neighbors that will never be forced.   If $S \backslash \{x\}$ is a zero forcing set of $T$ for some $x \in S'$, then $S' \backslash \{x\}$ is a zero forcing set of $T'$, contradicting minimality of $S'$.

Next, assume that Statement \ref{case:k} holds.  That is, suppose $S'$ is a minimal zero forcing set of $T'$ such that $ u \in S'$ for some $u \in V(T')$ with $uv \in E(T')$, and let $S=\left( S' \backslash \{u\} \right)  \cup L$. We see that $S$ is a zero forcing set of $T$, since without loss of generality, $l_1$ forces $v$, and since $S'$ restricted to any component of $T'$ not containing $u$ is a minimal zero forcing set of that component, $v$ forces $u$, after which $S'\cap V(T')$ is blue, and all of $T'$ will be forced.  For minimality of $S$, note that if we remove any vertex of $L$ from $S$, then $v$ cannot force $u$, and by minimality of $S'$, $T'$ will not be forced.  If $S \backslash \{x\}$ is a zero forcing set of $T$ for some $x \in S \backslash \{u\}$, then $S' \backslash \{x\}$ is a zero forcing set of $T'$, contradicting minimality of $S'$ and completing the reverse direction.

For the forward direction, first recall that if $S$ is a minimal zero forcing set of $T$, then $S$ contains $L$ or all but one vertex of $L$.  

First, suppose $l_j \notin S$ for some $j$, so $S=S'\cup L - \{l_j\}$ is a minimal zero forcing set of $T$.  We can see that $S'$ must be a zero forcing set on $T'$, since $v$ must force $l_j$ under $S$ and cannot force any vertices in $T'$.  If $S'$ is not minimal then there exists $x \in V(T')$ such that $S'\backslash\{x\}$ is a zero forcing set of $T'$.  But then $S \backslash \{x\}$ is a zero forcing set of $T$ as well, a contradiction of our assumption that $S$ is minimal, completing the case that $l_j \notin S$ for some $j$.  

Otherwise,  $S$ is a minimal zero forcing set of $T$ that contains all vertices in $L$.  In any realization of forcing chains under $S$, all but one leaf neighbor of $v$  must be a singleton forcing chain.  There exists a realization of forcing chains such that  one chain starts with the remaining leaf neighbor, goes through $v$ and 
continues to (without loss of generality) $u\in T_1$ where $T_1, T_2, \ldots T_m$ are the $m$ components of $T'$, and $m \geq 1$. That is, without loss of generality $l_1vu...$ is the beginning of a forcing chain.   Note that as before, if $m >1$ then for any $i$ with $1 < i \leq m$,  
$S\cap T_i$ is a minimal zero forcing set of $T_i$. By Lemma \ref{lem:pathcovers} the set $S-l_1 \cup \{v,u \}$ is a zero forcing set of $T$. We must be able to  delete $v$ since it is irrelevant by Proposition \ref{prop:pendent_not_minimal}. Hence  $S^*= S-l_1 \cup \{u \}$ is a zero forcing set of $T$.  To see that $S^*$ is  minimal, note  that if we can remove a vertex $x \neq u$ from $\left(S \cap T_1\right) \cup \{u\}$ and still produce a zero forcing set of $T_1$, then $S \backslash \{x\}$ is a zero forcing set of $T$, which contradicts the minimality of $S$. We can apply the same argument to $S^*$ that we used for the case that $l_j \notin S$  to complete the proof.
\end{proof}

\begin{cor}Let $T$ be a tree with vertex $v$ that has leaf neighbors $\{l_1, \dots , l_k\}$, where $k\geq2$. Define $T' = T-\{v, l_1, \dots , l_k\}$.  Then vertex $u \in V(T')$ is irrelevant in $T'$ if and only if it is irrelevant in $T$. \label{cor:gratisstatus}
\end{cor}

\begin{proof}
Suppose $u \in S \cap V(T')$ for a minimal zero forcing set $S$ of $T$.  Then by Theorem  \ref{thm:starremovalstillminimal2}, $u \in S'$ for some minimal zero forcing set $S'$ of $T'$, completing the forward direction.  

Suppose $u \in S'$ for some minimal zero forcing set of $T'$.  By Theorem \ref{thm:starremovalstillminimal2}, $S=S' \cup L \backslash \{l_1\}$ is a minimal zero forcing set of $T$, and $u\in S$ completing the proof.
\end{proof}

Suppose $G$ is a graph with a nonempty set $V_B$ of vertices with double pendants.  We call $G'$ a \emph{star reduction of $G$} if $G'$ is obtained by removing $V_B$ along with all leaf neighbors of $V_B$.    For example, in Figure \ref{figure:bvertices}, the star reduction of $T$ is the subforest with two components obtained by deleting the three darkest green vertices and their seven pendent neighbors from $T$.  

\begin{lem}
Under any minimal zero forcing set $S$ in a graph $G$, there is a realization of forcing chains so that each $B$-vertex is forced by one of its pseudoleaves. \label{lem:gencherrybelow}
\end{lem} 

\begin{proof}
Let $S$ be a minimal zero forcing set of $G$.  For any nonnegative integer $k$, recall that $B_k$ is the set of $B$-vertices that have double pendants after $k$ star reductions.  

We proceed by induction, inducting on the number of star reductions that are performed.  For the basis case, $k=0$ star reductions, we already know from  Proposition \ref{prop:pendent_not_minimal} 
that we have a realization of forcing chains with all vertices in $B_0$ forced by  their pseudoleaves.         

Assume that there is a realization of forcing chains of $S$ such that each vertex in $B_i$ is forced by one of its pseudoleaves for all $i$, $0 \leq i \leq k$.

Consider $B_{k+1}$, the nonempty set of $B$-vertices that have double pendants after $k+1$ star reductions.  Let $L \subseteq N(B_{k+1})$ be the set of pseudoleaves of $B_{k+1}$.  In $G$, the only neighbors of each vertex $u \in L$ are its lone neighbor in $B_{k+1}$ and a set of $B$-vertices that are removed in the first $k+1$ star reductions.  Let $\mathcal{B} = N(L) \backslash \{B_{k+1}\}$, that is the neighbors of $L$ in $B_i$ for $i \leq k$.  By the induction hypothesis, there is a realization of forcing chains such that each vertex in $\mathcal{B}$ is forced by one of its pseudoleaves.

Under $S$, each vertex $v \in B_{k+1}$ can force at most one of its neighbors.  Hence, for each $v \in B_{k+1}$ there exists a pseudoleaf $u \in L$ of $v$ such that $v$ cannot force $u$.   For $u$ to be forced, either $u \in S$, or $u$ is forced by  some vertex in $\mathcal{B}$.  By the argument above about $\mathcal{B}$, the only white neighbor of $u$ is $v$; hence $u$ can force $v$.  That is, there is a realization of the forcing chains such that each vertex in $B_i$ where $i \leq k+1$ is forced by its pseudoleaf, completing the proof. 
\end{proof}

Noting that any vertex that can be forced by another under a minimal zero forcing set $S$ cannot itself be in $S$ without violating minimality, we have the following corollary of Lemma \ref{lem:gencherrybelow}.

\begin{cor}
In any graph $G$, all $B$-vertices are irrelevant vertices. That is, no $B$-vertex is in any minimal zero forcing set of $G$.  
\label{cor:rootsnotinminimal}
\end{cor}

Finally, this section ends with a characterization of irrelevant vertices in trees. Corollary \ref{cor:rootsnotinminimal} establishes the reverse direction for all graphs; we prove the forward direction restricted to trees to complete the characterization.


\begin{thm} In a tree $T$, the vertex $v$ is an irrelevant vertex if and only if it is a $B$-vertex.\label{thm:treeminimalvertex}  
\end{thm}

\begin{proof}
The reverse direction is simply Corollary \ref{cor:rootsnotinminimal} restricted to trees.  

To prove the forward direction, suppose $v$ is an irrelevant vertex of $T$.  By Theorem \ref{thm:nodoublependants}, we know that $T$ must have a double pendant.  If $v$ is the neighbor of the double pendant, then $v$ is a $B$-vertex and we're done.  Otherwise by Corollary \ref{cor:gratisstatus}, we can remove the neighbor of the double pendant along with all of its leaf neighbors to produce $T'$, and $v$ is also irrelevant in $T'$.  By Theorem \ref{thm:nodoublependants}, $T'$ must have a double pendant. If $v$ is the neighbor of the double pendant, then we're done. Otherwise, we can again apply Corollary \ref{cor:gratisstatus}.  

Since $T$ has a finite number of vertices, eventually, if we apply Theorem \ref{thm:nodoublependants} and Corollary \ref{cor:gratisstatus} repeatedly, we must have that $v$ is the neighbor of the double pendant, and therefore a $B$-vertex.
\end{proof}

One direction of our characterization holds for graphs in general, not just trees. Also, all examples of irrelevant vertices in \cite{BongCarlsonCurtis} are $B$-vertices.  We end the section with an open question.

\begin{open}
Does the characterization in Theorem \ref{thm:treeminimalvertex} hold for all graphs, or just for trees?  
\end{open}

\section{Well-forced graphs} \label{Sec:well trees}
We now turn our attention back to characterizing well-forced graphs. In \cite{BorisCarlson}, the authors investigated the maximum cardinality of a minimal zero forcing set.  While they used different terminology from ours and primarily focused on other aspects of minimal zero forcing sets, we describe their results that directly relate to this paper in our terminology here.  First, they showed   \cite[Corollary 9]{BorisCarlson} that if $G$ is a cycle, empty graph, star, wheel, or complete graph, then $G$ is well-forced. They also showed  \cite[Proposition 10]{BorisCarlson} that for integers $a,b \geq 3$ that $K_a \vee \overline{K_b}$ is well-forced as is $K_a \vee C_b$, where $G \vee H$ denotes the \emph{join} of two graphs, that is, the graph formed from $G$ and $H$ by adding an edge between every vertex in $G$ and every vertex in $H$. By providing a construction, they proved that there are infinitely many graphs $G$ such that $G$ is well-forced but  $G \vee K_1$ is not \cite[Theorem 11]{BorisCarlson}.  We add the following graph families to the list of well-forced graphs.

\begin{thm}\label{path-not-wzf}
The path $P_n$ is well-forced if and only if $n\leq 3$.
\end{thm}

\begin{proof}
It is well known that $\Z(P_n)=1$.  If $n\geq 4$,  any adjacent pair of internal vertices form a minimal zero forcing set of order two.  
\end{proof}

\begin{prop}
If $G$ is  one of the following graphs, then it is well-forced.  
\begin{enumerate}
    \item a complete multipartite graph with at least two partite sets
    \item a cycle on at least three vertices with an additional single pendant
\end{enumerate}
\end{prop}

\begin{proof}
First suppose $G$ is a complete multipartite graph.  If $S$ is a zero forcing set, then $S$ omits at most two vertices from $V(G)$, the two vertices must be in different partite sets, and any set $S$ with this property is a zero forcing set.  Hence any set $S$ that omits exactly two vertices that are in two different partite sets is a minimal and minimum zero forcing set. 

Now suppose $G$ is a cycle with an additional pendant, and let $v$ be the vertex on the cycle with pendent neighbor $u$.  To be a zero forcing set, $S$ must contain either two adjacent vertices on the cycle, or $u$ and one other neighbor of $v$. To be minimal, $S$ must contain no other vertices. Thus $G$ is well-forced.  
\end{proof}

Given a graph $G$, the \emph{corona} of $G$ with $K_1$ denoted $G \circ K_1$ is the graph obtained by taking one copy of $G$ and $|V(G)|$ isolated vertices, assigning each isolated vertex a unique neighbor in $G$, and joining each of these pairs by an edge.  
\begin{prop}
For any connected graph $G$ with $|V(G)| \geq 2$, $G \circ K_1$ is not well-forced.
\end{prop}
\begin{proof}
Let $V(G) = \{v_1, \ldots, v_n\}.$  Let the leaf vertices generated by the corona $\{v'_1, \ldots, v'_n\}$.
Note that $V(G)$ is a minimal zero forcing set, and that $\{v'_1, \ldots, v'_{n-1}\}$ is  a zero forcing set with smaller cardinality than $V(G)$.  Hence $G\circ K_1$ is not well-forced.
\end{proof}


\subsection{Subgraphs that imply  a graph is not well-forced}
In this section we describe subgraphs whose presence guarantees that a graph that contains them will fail to be well-forced.  While these results help our characterization of well-forced trees,  the subgraphs described  need not be in a tree.  Thus we  also obtain conditions that preclude any graph from being well-forced.

\begin{prop} 
If $G$ has a pendent path $P$ of length four or more, then $G$ is not well-forced.  \label{prop:pendentpath}
\end{prop}

\begin{proof}
Starting with the leaf,  call the first three vertices of $P$ in order  $y, x, v$, and let $S$ be a minimum zero forcing set of $G$.  We can assume without loss of generality that $y \in S$ by using reversals, and we have $V(P_y) \cap S = \{y\}$.  

Let $S'=S \cup \{v,x\} \backslash \{y\}$.  Note that $|S'| =|S|+1$ and that $S'$ is a zero forcing set.  We show that $S'$ is minimal.  
Suppose we delete $v$ or $x$ from $S'$.  Then $y$ will not be forced, because if it were, then we would have that $S \backslash\{y\}$ is a zero forcing set, a contradiction.  Also, if $S'\backslash \{z\}$ is a zero forcing set  for some $z \in S \backslash V(P_y)$, then $S \backslash \{z\}$ is a zero forcing set as well, a contradiction.  
Hence $S'$ is minimal, giving us that $G$ is not well-forced. 
\end{proof}

\begin{prop}
 Suppose $G$ has a pendent generalized star $R$ whose legs all have length at least two.  
Then $G$ is not well-forced. \label{prop:pgsnoshortlegs}
\end{prop}

\begin{proof}
Let $S$ be a minimum zero forcing set of $G$.  By Corollary \ref{cor:nocenter}, the center vertex $v$ of $R$ is not in $S$.  Then $v$ is forced by either by a vertex on $R$, or by the neighbor $u$ of $v$ that is not on $R$. First assume $u$ does not force $v$.  Then since $S$ is a minimum zero forcing set, a leaf on some leg of $R$ is in $S$, and all legs but at most one either have a leaf in $S$, or the vertex adjacent to $v$ in $S$.   Let $S'$ consist of all vertices in $S$ from $G \backslash R$, along with $v$ and the vertex adjacent to $v$ on each leg that has a vertex in $S$.  Then $S'$ is a minimal zero forcing set with cardinality $|S|+1$, and $G$ is not well-forced.

Next suppose $u$ forces $v$ under $S$. Then on each leg except one, either the leaf or the vertex adjacent to $v$ is in $S$. Since either choice will produce a minimum zero forcing set, assume without loss of generality that the vertex adjacent to $v$ on each of these legs is in $S$.  For the final leg $L$ without a vertex from $S$, we must have that $v$ forces its neighbor on $L$, and the forcing process continues to the leaf of  $L$.  Consider the reversal $\tilde{S}$ of $S$.  Then $\tilde{S}$ contains all leaves of $T$, so we can construct $\tilde{S}'$ as we constructed $S'$ above, completing the proof that $G$ is not well-forced.   
\end{proof}

\begin{thm}
Let $G$ be a generalized star with  center vertex $u$ of degree $k\geq 2$. Then $G$ is well-forced if and only if $G$ has at least two legs of length one and all legs of length at most three.  \label{thm:genstarwellchar}
\end{thm}

\begin{proof}
Recall that $\Z(G)=k-1$.  

For the forward direction, note that if $G$ has at most one leg of length one, we can take the center vertex $u$ and $k-1$ of its degree-two neighbors to form a minimal zero forcing set $S$ with $|S| = k$.  If $G$ has a leg of length four or more, then $G$ is not well-forced by Proposition \ref{prop:pendentpath}, completing the forward direction.

For the reverse direction, suppose $G$ is a generalized star with at least two legs of length one and all legs of length at most three.  Let $S$ be a minimal zero forcing set of $G$. Note that there is at least one leg $L$ of length one whose sole vertex lies in $S$.  This means that the center vertex $u$ can be forced by this vertex. Thus, $u \notin S$.
By Observation \ref{obs:allbutoneleg}, there exist $k-1$ legs $L$ such that $L \cap S \neq \emptyset$.  So $|S| \geq k-1.$    Let $L$ be an arbitrary leg of $G$ with $L \cap S \neq \emptyset$.  If $L$ is length one, it is clear that $|L \cap S| \leq 1$.    If $L$ has length two or three, then either end vertex of $L$ on its own is sufficient to force the remaining vertices in $L$  (since $u$ will be forced by a length-one leg).  Thus, $|L \cap S| \leq 1$.  Hence, $|S|=\Z(G)=k-1$, completing the proof.
\end{proof}

\subsection{A characterization and algorithm for well-forced trees}

In this section, we describe a process that preserves the property of being well-forced, which allows us to characterize  well-forced trees.  
For any $v \in V(T)$, let $L(v)$ denote $N(v) \cap \{u \in V(T) : \deg(u) = 1\}$.  Let $L[v]=L(v)\cup\{v\}$. 

\begin{prop}
Let $T$ be a tree with  $|L(v)|\geq 2$ for some $v \in V(T)$.
Then $T$ is well-forced if and only if each component of the graph induced by $V(T) \backslash L[v]$ is well-forced. \label{prop:starremoval}
\end{prop}
\begin{proof}
For the forward direction, we prove the contrapositive. 

Let the components of the graph induced by $V(T) \backslash L[v]$ be $T_1, T_2, \ldots T_k, k\geq 1$.  Suppose that $T_1$ is not well-forced.  Let $S_i$ be a minimum zero forcing set of $T_i$ for each $i$, $1 \leq i \leq k$. Since $T_1$ is not well-forced, we can find a minimal zero forcing set $S'_1$ such that $S'_1$ is a minimal zero forcing set of $T_1$ and $|S_1| < |S'_1|$.

Let $L(v)=\{\ell_1, \ell_2, \ldots, \ell_j\}$ where $j \geq 2$. Let $S=S_1 \cup S_2 \cup \cdots \cup S_k \cup L(v) \backslash \{\ell_j\}$ and $S'=\left( S \backslash S_1 \right) \cup S'_1$.  Then by Observation \ref{obs:disconnected} and Theorem \ref{thm:starremovalstillminimal2}, $S, S'$ are minimal zero forcing sets of $T$ with $|S| < |S'|$, completing the forward direction.



For the reverse direction, suppose each component $T_i$ of the graph $T'$ induced by $V(T) \backslash L[v]$ is well-forced.  Let $S$ be a minimal zero forcing set of $T$.  Then by Theorem \ref{thm:starremovalstillminimal2}, either $S$ contains all vertices in $L(v)$ and all but one vertex of a minimal zero forcing set of $T'$, or $S$ contains all but one vertex in $L(v)$, and all vertices of a minimal zero forcing set of $T'$.  By Observation \ref{obs:disconnected} and the assumption that all components of $T'$ are well-forced, we find that $|S|$ is fixed and $T$ is well-forced.

\end{proof}

\begin{thm}
If $T$ is a tree with $|V(T)|>2$ and no double pendants, then $T$ is not well-forced.   \label{thm:nostars}
\end{thm}

\begin{proof}
If $T$ is a path or a generalized star, then the result holds by  Theorem \ref{path-not-wzf} or \ref{thm:genstarwellchar} respectively.

Otherwise, $T$ has one or more pendent generalized stars  \cite{fallat2007minimum}, each with at most one leg of length one.  If any pendent generalized star has  a leg of length four or more or no legs of length one,  then by Propositions \ref{prop:pendentpath} and \ref{prop:pgsnoshortlegs}, $T$ is not well-forced. 

Hence, the remaining case is that each pendent generalized star on $T$ has one leg of length one, and all other legs of length two or three.  
Pick any pendent generalized star and call it $R$.  Let its center vertex be called $c$, and let $T'=T \backslash R$.  We know that $c$ is in some  minimal zero forcing set of $T$ by Theorem \ref{thm:nodoublependants}.  Let $S$ be a minimal zero forcing set with $c \in S$.  By Lemma \ref{centerminimalhasneighbors}, $S\cap R$ consists of $c$ and at least $\ell-1$ neighbors of $c$ on $R$ where $\ell$ is the number of legs of $R$.  

Let $S'=\left(S\cap T'\right) \cup L$ where $L$ is the set of leaves on legs of $R$ that have a vertex in $S$. Then $S'$ is also a zero forcing set with $|S'|=|S|-1$, and hence $T$ is not well-forced.
\end{proof}

Finally, we present our characterization of well-forced trees.  We perform a \emph{star removal} on $T$ by  identifying a vertex $v$ with $|L(v)|\geq 2$, and creating a new tree $T \backslash L[v]$.  

\begin{thm}
A tree $T$ is well-forced if any only if performing star removals on $T$ until no more star removals can be performed results in a possibly empty set of copies of $K_2$.  \label{thm:wellcharacterization}
\end{thm}

\begin{proof}
Suppose $T$ is a well-forced tree.  Then $T$ contains a vertex $v$ with a double pendant by Theorem \ref{thm:nostars}. Let $T_1= T\backslash L[v]$. By Proposition \ref{prop:starremoval}, $T_1$ is well-forced. Repeat this process on $T_1$ and the graphs that follow until no more star removals can be done. Call the resulting forest $F$. By Proposition \ref{prop:starremoval} $F$ is well-forced, and no vertex has a double pendant (otherwise we would perform another star removal). By Theorem \ref{thm:nostars}, either $|V(F)|=0$ or each component has at most two vertices. Note that no component can be a singleton vertex $u$, since this implies the neighbors of $u$ were removed, and must have been vertices with double pendants, but then $u$ would have been a pendant after some  of star removals, and therefore removed as well. Thus if $|V(F)|\neq0$, each component of $F$  is isomorphic to $K_2$.

Now suppose that $T$ is a tree with the property that after performing star removals until no more can be performed, all components of the resulting graph $F$ are copies of $K_2$. Any minimal zero forcing set of $K_2$ must contain exactly one vertex. Hence, each component of $F$ is well-forced. Therefore by Proposition \ref{prop:starremoval}, the original tree $T$ is also well-forced.  
\end{proof}

Theorem \ref{thm:wellcharacterization} leads to the following algorithm for determining whether or not a tree is well-forced.  

\begin{alg}
Well-forced tree
\end{alg}
\noindent {\bf Input:} A tree $T$\\
{\bf Output:} True/False variable \emph{WELL}
\begin{enumerate}
\item While there exists $v \in V(T)$ with $|L(v)| \geq 2$:
\begin{enumerate}[A.]
\item Choose a vertex $v_0$ with $|L(v_0)| \geq 2$.
\item Let $T=T \backslash L[v_0]$.     
\end{enumerate}
\item If $T= \emptyset$ or $T$ is the union of copies of $K_2$:\\
 Then \emph{WELL}=True;\\
 Else \emph{WELL}=False.
\end{enumerate}

\end{document}